\documentclass[onecolumn, 11pt]{article}
\usepackage[utf8]{inputenc}
\usepackage[left = 2cm, right=2cm, top=2cm]{geometry}
\usepackage{setspace}
\usepackage{comment}

\usepackage{amsmath, amssymb}
\usepackage{graphicx}
\usepackage{dcolumn}
\usepackage{bbm}
\usepackage{hyperref}
\usepackage{csquotes}
\usepackage{float}
\usepackage{caption}
\usepackage{subcaption}
\usepackage{amsthm}
\usepackage{tabularx}
\usepackage{color}

\newtheorem{theorem}{Theorem}
\newtheorem{lemma}{Lemma}

\newcommand{\E}{\mathbb{E}}
\newcommand{\prob}{\mathbb{P}}
\newcommand{\Reals}{\mathbb{R}}
\newcommand{\Nats}{\mathbb{N}}

\newcommand{\niso}{N_{\textrm{iso}}}
\newcommand{\nisotilde}{\widetilde{N}_{\textrm{iso}}}

\newcommand{\iso}{W_{m, L}}
\newcommand{\PPP}{\mathcal{P}}
\newcommand{\htilde}{\tilde{h}}
\newcommand{\Htilde}{\widetilde{H}}
\newcommand{\Gtilde}{\mathcal{G}_{\tilde{h}^L}(\mathcal{P}_L)}

\title{The Distribution of the Number of Isolated Nodes in the 1-Dimensional Soft Random Geometric Graph}
\author{Michael Wilsher
\and Carl P. Dettmann \and Ayalvadi J. Ganesh
}
\date{\textit{School of Mathematics, University of Bristol, Bristol BS8 1UG, UK}}

\begin{document}
\maketitle

\begin{abstract}
We study the number of isolated nodes in a soft random geometric graph whose vertices constitute a Poisson process on the torus of length $L$ (the line segment $[0, L]$ with periodic boundary conditions), and where an edge is present between two nodes with a probability which depends on the distance between them. Edges between distinct pairs of nodes are mutually independent. In a suitable scaling regime, we show that the number of isolated nodes converges in total variation to a Poisson random variable. The result implies an upper bound on the probability that the random graph is connected.
    
\end{abstract}


\section{Introduction}

Random geometric graphs (RGGs) were introduced in \cite{gilbert1959random} as a model for communication networks with limited connection range and have subsequently been widely used to model networks with a spatial element; see, e.g., \cite{penrose2003} and references therein. In Gilbert's model, nodes or vertices are randomly distributed within some space, typically $\Reals^d$ or a bounded subset of it, and an edge is placed between nodes if their mutual distance is smaller than a specified threshold. We shall refer to this model as the hard RGG, as the connection probability has a hard cutoff as a function of the distance between nodes. 
A generalisation is to allow the edge between a pair of nodes separated by distance $r$ to be present with some probability $H(r)$, independent of all other edges. This model has been termed the Soft RGG \cite{penrose2016connectivity}, the random connection model (RCM) \cite{meester1996continuum}, or a Waxman graph \cite{waxman1988routing}. In this work we refer to them as Soft RGGs and term $H(\cdot)$ the connection function. The hard RGG is a special case, obtained by setting $H(r) = 1$ for $r \leq r_c$ and $H(r) = 0$ otherwise; $r_c > 0$ is a parameter of the model.

The one-dimensional version of this model is motivated by vehicular ad-hoc networks (VANETs), which are expected to be essential for autonomous vehicles; these will be fitted with on-board radios to transmit information such as location, velocity, and hazard warnings between vehicles. The road is represented by a line, the nodes represent vehicles and an edge between two nodes indicates that the two vehicles can communicate directly. One key question is connectivity: When is every vehicle in the network able to communicate with every other, either via a single- or multi-hop path? A necessary condition for full connectivity is the absence of isolated nodes, namely nodes that are not connected to any other in the graph. In the 2-D and 3-D versions of the soft RGG model, it has been shown in a suitable asymptotic regime that the soft RGG is connected if and only if there are no isolated nodes \cite{mao2012connectivity, penrose2016connectivity}. It was further shown in \cite{mao2012connectivity} that the the number of isolated nodes in a 2-D soft RGG can be well approximated by a Poisson distribution. The result was extended in \cite{penrose2016connectivity} to dimension three and greater; also see \cite{ganesh2007} for an analogous result for small-world networks. There has been little work to date on the 1-D model, which differs in important respects from those in two or more dimensions. In particular, the dominant reason for disconnection in 1-D hard RGGs is the presence of uncrossed gaps rather than of isolated nodes. Soft RGGs in 1-D were studied in~\cite{wilsher2020connectivity}, where it was shown that isolated nodes dominate uncrossed gaps as a cause of disconnection. The threshold for the emergence of isolated nodes was established in~\cite{wilsher2020connectivity}, and a Poisson limit was conjectured for the number of isolated nodes at the threshold (Conjecture 4.1). We prove this conjecture here.



We now describe the model studied in this paper and our assumptions. We consider a sequence of networks indexed by a parameter $L\in \Reals_+$, which tends to infinity along the sequence. The nodes or vertices of the network constitute a Poisson point process (PPP) of unit intensity on $[0,L]$, which we denote $\mathcal{P}_L$. If two nodes are located a distance $r$ apart, the edge between them is present with probability $H^L(r)= H(r/R_L)$, independent of all other edges; here, $H:\Reals^+\to [0,1]$ is a given connection function, and $R_L>0$ is a scaling parameter to be specified.

In order to avoid inessential technicalities associated with the boundaries, we shall study the model with periodic boundary conditions.  In other words, we turn $[0,L]$ into a circle or torus by identifying $0$ with $L$. We denote the circular distance on $[0,L]$ by $\rho^L$, i.e., $\rho^L(x,y)=\min \{ |x-y|, L-|x-y| \}$. Finally, we define the connection probabilities to be $h^L(x,y) = H^L(\rho^L(x,y))$. We denote by $\mathcal{G}_{h^L}(\mathcal{P}_L)$ the graph with vertex set $\mathcal{P}$, and independent edges generated with connection probabilities $h^L$. This is the graph model we study. We believe that the same results hold if $[0,L]$ is treated as a line segment rather than a circle, and that this can be established by analysing boundary effects separately, as in~\cite{dettmann2016random}.

We make the following assumptions about the connection function.

\noindent \textbf{Assumptions:} Let $\| H\|_1$ and $\| H\|_2$ denote the $L^1$ and $L^2$ norms of $H:\Reals_+ \to [0,1]$, i.e., 
$$
\| H\|_1=\int_0^{\infty} H(x)dx, \quad \| H\|_2^2 = \int_0^{\infty} H(x)^2 dx.
$$
We assume that $\| H\|_1<\infty$ and $\| H\|_2^2< \| H\|_1$.

The first assumption says that $H$ is integrable, and is required for the mean number of neighbours of a node to remain bounded as $L$ tends to infinity. We expect this to hold in real-world networks. Next, observe from the definition that $\| H\|_2^2 \leq \| H\|_1$, since $H(x) \in [0,1]$ for all $x$ as it is a probability. Hence, the second assumption, which asserts that this inequality is strict, is a mild one. It is satisfied whenever the set, $\{ x\in \Reals_+: 0<H(x)<1 \}$, where the connection probability is strictly between 0 and 1, has positive Lebesgue measure. Nevertheless, the assumption excludes the connection function of a hard RGG, which is $\{ 0,1\}$-valued.

We now state our main result, which resolves a conjecture in \cite{wilsher2020connectivity}.

\begin{theorem}
\label{thm:IsoNodesPoisson_torus}
Fix $\tau \in \Reals_+$ and take $R_L=\ln(\tau L)/2\| H\|_1$. Let $\niso$ denote the number of isolated nodes in the soft RGG $\mathcal{G}_{h^L} (\mathcal{P}_L)$. Its dependence on $L$ and  $\tau$ has been suppressed for notational convenience. As $L$ tends to infinity, the random variable $\niso$ converges in total variation distance to a Poisson random variable with mean $1/\tau$. In particular, $\prob(\niso=0)$ tends to $e^{-1/\tau}$.
\end{theorem}

As the soft RGG is disconnected if there are any isolated nodes, an immediate corollary of the theorem is that the probability that $\mathcal{G}_{h^L}(\mathcal{P}_L)$ is connected is asymptotically bounded above by $e^{-1/\tau}$, in the scaling regime considered in the theorem. This upper bound would be tight if isolated nodes were the dominant cause of disconnection in this random graph model, as conjectured in \cite{wilsher2020connectivity} under the mild additional condition that the connection function has unbounded support. Resolving this conjecture remains an open problem, as does extending the analysis of this paper to point process models other than the Poisson process.

\section{Proofs} 
\label{sec:toroidal_conn_fun}

Denote the total variation distance between probability distributions $\mu$ and $\nu$ on $\Reals$ by $d_{TV}(\mu,\nu)$. With a slight abuse of notation, we shall write $d_{TV}(X,Y)$ for random variables $X$ and $Y$ to denote the total variation distance between their probability distributions. The following bound on the total variation distance between random variables $X$ and $Y$ defined on the same probability space is well-known and elementary:
\begin{equation}
    \label{eq:vardist_coupling_bd}
    d_{TV}(X,Y) \leq \prob(X\neq Y).
\end{equation}

The proof of Theorem \ref{thm:IsoNodesPoisson_torus} proceeds through a sequence of lemmas. Our first result approximates the number of isolated nodes in our model with one in which the connection function is truncated. This step is not needed if the connection function has bounded support.


\begin{lemma}
\label{lem:niso_nisotilde_lim}
Fix $\alpha>0$, and define the truncated connection function $\Tilde{h}^L(x,y) = h^L(x, y)\mathbbm{1}\lbrace \rho^L(x, y) \leq R_L^{1+1/\alpha} \rbrace.$
Denote by $\niso$ and $\nisotilde$ the number of isolated nodes in $\mathcal{G}_{h^L}(\mathcal{P}_L)$ and $\Gtilde$ respectively. Then,
$$
\lim_{L\rightarrow\infty}d_{TV}(\niso, \nisotilde) = 0.
$$
\end{lemma}

\begin{proof}
We can couple $\mathcal{G}_{h^L}(\mathcal{P}_L)$ and $\Gtilde$ by first generating $\mathcal{G}_{h^L}(\mathcal{P}_L)$, and then removing any edges of length at least $R_L^{1 + 1/\alpha}$. Observe that $\nisotilde \geq \niso$, since removing edges cannot reduce the number of isolated nodes. Therefore, it follows from Markov's inequality that
\begin{equation}
\label{eqn:markov_to_show_dtv}
\prob(|\niso - \nisotilde| \geq 1) \leq \E[|\niso - \nisotilde|] = \E[\nisotilde] - \E[\niso].
\end{equation}
Using the expression for the expected number of isolated nodes in a soft RGG derived in \cite[eqn. (9)]{wilsher2020connectivity}, we get
$$
    \E[\niso] = L\exp\left(-2R_L\int_0^{L/2R_L} H(r)dr\right) 
    = L\exp\left(-\frac{\ln(\tau L)}{\|H\|_1}\int_0^{L/2R_L} H(r)dr \right). 
$$
Now, $L/2R_L$ tends to infinity as $L$ tends to infinity, and so 
$\int_0^{L/2R_L} H(r)dr$ tends to $\| H\|_1=\int_0^{\infty} H(r)dr$. Hence, it follows from the above that 
\begin{equation}
\label{eqn:niso_inf}
    \E[\niso] 
    \xrightarrow{L\rightarrow\infty} \frac{1}{\tau}.
\end{equation}
Similarly, since $R_L=\frac{\ln(\tau L)}{2\| H\|_1}$ tends to infinity with $L$, we have 
\begin{equation}
\label{eqn:nisotilde_inf}
\begin{aligned}
    \E[\nisotilde] &= L\exp\left(-2\int_0^{L/2} \tilde{h}^L(0,r) dr\right)  = L\exp \left( -2\int_0^{R_L^{1+1/\alpha}} H\Bigl( \frac{r}{R_L} \Bigr)dr \right) \\
    & = L\exp\left(-\frac{\ln(\tau L)}{\|H\|_1}\int_0^{R_L^{1/\alpha}} H(r)dr\right) \xrightarrow{L\rightarrow\infty} \frac{1}{\tau}.
\end{aligned}
\end{equation}
It follows from \eqref{eqn:niso_inf} and \eqref{eqn:nisotilde_inf} that $\E[\nisotilde]-\E[\niso]$ tends to zero as $L$ tends to infinity. The claim of the lemma follows from \eqref{eq:vardist_coupling_bd} and \eqref{eqn:markov_to_show_dtv}.
\end{proof}

Henceforth, we shall work with $\Gtilde$, the soft RGG with truncated connection function. We shall use the Chen-Stein method for Poisson approximation described in \cite{barbour1992poisson}; as it requires a discrete index set, we fix $m\in \Nats$ and discretise the torus $[0,L]$. into $mL$ segments of width $1/m$. (Assume for convenience that $L$ is an integer. Otherwise, we need one segment of a different width, which does not fundamentally alter the analysis.) Denote the $i$th segment by $A_{i}$, where $i$ takes values in the index set $\Gamma = \lbrace 1, 2, ..., mL \rbrace$. Denote by $\mathcal{P}_L$ and $\mathcal{I}$ the vertex set and the set of isolated nodes of the random graph $\Gtilde$, and by $\mathcal{P}_L(A_i)$ and $\mathcal{I}(A_i)$ the subsets of these nodes lying within $A_i$. Define
\begin{equation}
    \label{eqn:isol_defs}
    \begin{aligned}
    J_{i} &:= \mathbbm{1}\left\lbrace |\mathcal{P}_L(A_{i})| = 1\right\rbrace, \\
    I_{i} &:= \mathbbm{1}\left\lbrace |\mathcal{P}_L(A_{i})| = 1 \textrm{ and } |\mathcal{I}(A_i)|=1 \right\rbrace, \\
    W_{m, L} &:= \sum_{i \in \Gamma}I_i, 
    \end{aligned}
\end{equation}
for $i \in \Gamma$. Denote the centre of the segment $A_{i}$ by $x_{i}$. (Although $A_i, I_i, J_i,$ and $x_i$ all implicitly depend on $m$, this dependence is suppressed for notational convenience.)
Our next result states that the number of isolated nodes in well-approximated by $W_{m,L}$ when $m$ is large.
\begin{lemma}
\label{lem:discrete_approx}
Let $I_i$ denote the indicator that the $i^{\rm th}$ segment of $[0,L]$ contains exactly one node, and that node is isolated in $\Gtilde$. Let $W_{m,L}$ denote the sum of these indicators, as defined in \eqref{eqn:isol_defs}. Let $\nisotilde$ denote the total number of isolated nodes in $\Gtilde$. Then, for any fixed $L>0$,
$$
d_{TV}(\nisotilde, W_{m,L}) \to 0 \mbox{ as } m\to \infty.
$$
\end{lemma}
\begin{proof}
Notice that the random variable $W_{m,L}$ is exactly the same as the number of isolated nodes unless there is a segment containing two or more nodes. Now, the number of nodes in a segment of length $1/m$ has a Poisson distribution with mean $1/m$. Hence, it follows from \eqref{eq:vardist_coupling_bd} and the union bound that 
\begin{equation*}
    \begin{aligned}
        d_{TV}(\nisotilde,W_{m,L}) &\leq \prob(\exists i\in \Gamma: |\mathcal{P}_L(A_i)|\geq 2) \\
        &\leq \sum_{i\in \Gamma} \Bigl( 1 -e^{-1/m} -\frac{1}{m} e^{-1/m}\Bigr) \\
        &= mL \Bigl( 1 -e^{-1/m} -\frac{1}{m} e^{-1/m}\Bigr),
    \end{aligned}
\end{equation*}
which tends to zero as $m$ tends to infinity, as claimed. 
\end{proof}

In light of the above lemma, it suffices to establish a Poisson approximation for the random variables $W_{m,L}$. To this end, we define the ``neighbourhood" of an index $i\in \Gamma$ by 
\begin{equation} \label{eqn:Bi_def}
B_i = \left\lbrace j \in \Gamma: \rho^L(x_i, x_j) \leq 3R_L^{1 + 1 / \alpha} \right\rbrace.
\end{equation}
We also define the quantities
\begin{equation}
\label{eqn:def_b123}
    \begin{aligned}
        &p_i := \E[I_i], \quad p_{ij} := \E[I_i I_j], \\
        &b_1 := \sum_{i \in \Gamma} \sum_{j \in B_i} p_i p_j, \quad b_2 := \sum_{i \in \Gamma} \sum_{j \in B_i \backslash \lbrace i \rbrace} p_{ij}, \quad b_3 &:= \sum_{i \in \Gamma}\E \left[| \E\left[I_i | (I_j : j \in \Gamma\backslash B_i) \right] - p_i |\right].
    \end{aligned}
\end{equation}
We shall use of the following result on Chen-Stein approximation.
\begin{lemma}(\cite[Theorem 1.A]{barbour1992poisson}) 
\label{lem:Chen-Stein}
Let $Po(\lambda)$ denote a Poisson random variable with mean $\lambda$. Let $W_{m,L}$ be defined as in (\ref{eqn:isol_defs}) and $b_1$, $b_2$ and $b_3$ as in (\ref{eqn:def_b123}). Then,
$$
d_{TV}\left(\iso, Po\left(\E[\iso]\right) \right) \leq \min\left(1, \frac{1}{\E[\iso]}\right) (b_1 + b_2 + b_3).
$$
\end{lemma}
Thus, in order to establish a Poisson approximation for $\iso$, we need to bound the terms $b_1,b_2,b_3$. 

\begin{lemma} \label{lem:b3_bound}
Let $B_i$ be defined as in \eqref{eqn:Bi_def} and $b_1,b_2,b_3$ as in \eqref{eqn:def_b123}. Then, $b_3=0$ for all $m$ sufficiently large.
\end{lemma}


\begin{proof}
The connection function $\tilde h_L$ is defined by truncating $h_L$ at $R_L^{1+1/\alpha}$. Hence, the event for which $I_i$ is the indicator depends only on nodes within distance $R_L^{1+1/\alpha}$ of the segment $A_i$. By the same reasoning, the isolation or otherwise of nodes at distance greater than $2R_L^{1+1/\alpha}$ from this segment is independent of the Poisson point process in the set $\{ x:\min_{y\in A_i}\rho^L(x,y) \leq R_L^{1+1/\alpha} \}$. Hence, $I_i$ is jointly independent of $I_j$ for all $j$ such that $\min_{x\in A_i,y\in A_j} \rho^L(x,y) > 2R_L^{1+1/\alpha}$. In particular, if $m$ is large enough that $1\frac{1}{m} < R_L^{1+1/\alpha}$, then $I_i$ is independent of the collection of random variables, $\{ I_j, j\in \Gamma\setminus B_i \}$. Hence, $b_3=0$ for all $m$ sufficiently large.
\end{proof}

\begin{lemma} \label{lem:b1_bound}
Let $B_i$ be defined as in \eqref{eqn:Bi_def} and $b_1,b_2,b_3$ as in \eqref{eqn:def_b123}. Then, $b_1$ tends to zero as we let $m$ tend to infinity, followed by $L$.
\end{lemma}

\begin{proof}
Observe from the definition of $I_i$ that 
\begin{equation} \label{eqn:pi_def}
\begin{aligned}
p_i &= \prob(J_i=1) \prob(I_i=1|J_i=1) \\
&= \frac{1}{m} e^{-1/m}\prob(v\in \mathcal{P}_L(A_i) \mbox{ is isolated in $\Gtilde$}| J_i=1 ),
\end{aligned}
\end{equation}
where $v\in \mathcal{P}_L(A_i)$ exists and is unique since $J_i=1$. Say $v$ is located at $(i/m)+x$, where $0\leq x<1/m$. Now, the set of nodes to which $v$ is connected form an inhomogenous Poisson process on $[0,L]\setminus A_i$, of intensity $\tilde h^L(i/m+x,y)$; $v$ is isolated only if this set is empty. Thus, by translation invariance of the connection function, the probability that $v$ is isolated is given by 
\begin{align*}
    &\prob(v\mbox{ is isolated}) = e^{-\kappa_x},\mbox{ where } \\ 
    &\kappa_x = \int_{[0,L]\setminus [0,1/m]} \tilde h^L(x,y)dy 
    \xrightarrow{m\to \infty} \int_0^L \tilde h^L(x,y) dy,
\end{align*}
since $\tilde h^L$ is bounded above by 1. We have not made the dependence of $\kappa_x$ on $m$ and $L$ explicit in the notation. Suppose $L$ is large enough that $L/2R_L \geq R_L^{1+1/\alpha}$. Then, using the translation invariance of $\tilde h^L$ once more, we have 
$$
\lim_{m\to \infty} \kappa_x = \int_0^L \tilde h^L(0,y)dy \\
= 2\int_0^{R_L^{1+1/\alpha}} H \Bigl( \frac{y}{R_L} \Bigr) dy \\ 
= 2R_L \int_0^{R_L^{1/\alpha}} H(y)dy,
$$
which does not depend on $x$. Substituting in \eqref{eqn:pi_def}, we get that 
$$
\lim_{m\to \infty} mp_i = \lim_{m\to \infty} e^{-1/m}e^{-\kappa_x} = e^{-2R_L \int_0^{R_L^{1/\alpha}} H(y)dy.}
$$
Substituting in the definition of $b_1$, and noting that $p_i$ does not depend on $i$, we see that 
\begin{align*}
    \lim_{m\to \infty} b_1 
    &= \lim_{m\to \infty} |\Gamma| |B_i| p_i^2 
    = \lim_{m\to \infty} 6m^2 LR_L^{1+1/\alpha} p_i^2 \\ 
    &= 6L \Bigl( \frac{\ln(\tau L)}{2\| H\|_1} \Bigr)^{1+1/\alpha} 
    \exp \Bigl( -\frac{2\ln(\tau L)}{\| H\|_1} \int_0^{R_L^{1/\alpha}} H(y)dy \Bigr).
\end{align*}
Now, $\int_0^{R_L^{1/\alpha}} H(y)dy$ tends to $\| H\|_1$ as $L$, and hence $R_L$, tends to infinity. So, it follows that 
$$
\lim_{L\to \infty} \lim_{m\to \infty} b_1 
= \lim_{L\to \infty} \frac{6L}{(\tau L)^2} \Bigl( \frac{\ln(\tau L)}{2\| H\|_1} \Bigr)^{1+1/\alpha} = 0.
$$
This completes the proof of the lemma.
\end{proof}

\begin{lemma} \label{lem:b2_bound}
Let $B_i$ be defined as in \eqref{eqn:Bi_def} and $b_2$ as in \eqref{eqn:def_b123}. Then, $b_2$ tends to zero as we let $m$ tend to infinity, followed by $L$.
\end{lemma}

\begin{proof}
Consider two nodes located at $x,y\in [0,L]$. The set of all other nodes to which at least one of them has an edge constitutes an inhomogeneous PPP on $[0,L]$, with intensity $\phi(\cdot,\{ x,y\})$ given by 
\begin{equation} \label{eqn:phi_def}
\phi(z,\{ x,y\}) = 1-(1-\htilde^L(z,x))(1-\htilde^L(z,y)).
\end{equation}

Fix $i, j\in \Gamma$, $i\neq j$, and condition on the event that $J_i=1$ and $J_j=1$, i.e., that there is a unique point of the PPP, $\mathcal{P}_L$, on each of the segments $A_i$ and $A_j$. Denote their positions by $x$ and $y$. The set of nodes to which these might be connected, besides each other, constitutes a PPP on $[0,L]\setminus(A_i\cup A_j)$ with intensity $\phi$ defined above. Hence, the probability that both these nodes are isolated is given by 
$$
(1-\htilde^L(x,y))\exp \Bigl( -\int_{[0,L]\setminus (A_i\cup A_j)} \phi(z,\{ x,y\})dz \Bigr).
$$
Now, by well-known properties of Poisson point processes, the unique points of the homogeneous PPP on the segments $A_i$ and $A_j$ are independently and uniformly distributed within them. Hence, 
\begin{align*}
    p_{ij}
    &= \prob(J_i=1,J_j=1) \prob(I_i=1,I_j=1|J_i=1,J_j=1) \\ 
    &= \frac{1}{m^2}e^{-2/m} \int_{x\in A_i}\int_{y\in A_j} m^2 (1-\htilde^L(x,y)) \exp \Bigl(-\int_{[0,L]\setminus (A_i\cup A_j)} \phi_{\htilde^L}(z,\{ x,y\}) dz \Bigr) dxdy \\
    &\leq \int_{x\in A_i}\int_{y\in A_j} \exp \Bigl( -\int_{[0,L] \setminus (A_i\cup A_j)} \phi_{\htilde^L}(z,\{ x,y\}) dz \Bigr) dxdy
\end{align*}
Substituting in the definition of $b_2$, we obtain that 
$$
    b_2 \leq \int_{x\in [0,L]}\int_{y: \rho^L(x,y) \leq 3R_L^{1+1/\alpha}} \exp \Bigl(-\int_{[0,L]\setminus (A_i(x)\cup A_j(y))} \phi_{\htilde^L}(z,\{ x,y\}) dz \Bigr) dxdy,
$$
where we write $A_i(x)$ and $A_j(y)$ to denote the segments in which $x$ and $y$ lie; allowing them to lie in the same segment yields an upper bound, as does dropping the $(1-\tilde h^L(x,y))$ term. Now, $$
\int_{[0,L]\setminus (A_i(x)\cup A_j(y))} \phi_{\htilde^L}(z,\{ x,y\}) dz \xrightarrow{m\to \infty} \int_{[0,L]} \phi_{\htilde^L}(z,\{ x,y\}) dz,
$$
since $\phi_{\htilde^L}$ is bounded above by 1. Hence, we conclude using the translation invariance of $\htilde^L$ that 
\begin{equation} \label{eqn:b2_bd}
    \begin{aligned}
    \limsup_{m\to \infty} b_2 
    &\leq \int_{x\in [0,L]}\int_{y: \rho^L(x,y) \leq 3R_L^{1+1/\alpha}} \exp \Bigl( -\int_{[0,L]} \phi_{\htilde^L}(z,\{ x,y\}) dz \Bigr) dxdy \\
    &= 2L \int_{y=0}^{3R_L^{1+1/\alpha}} \exp \Bigl( -\int_0^L \phi_{\htilde^L}(z,\{ 0,y\}) dz \Bigr) dy. 
    \end{aligned}
\end{equation}
Substituting for $\phi_{\htilde^L}$ from \eqref{eqn:phi_def}, we have 
\begin{equation} \label{eqn:phi_integral_0}
    \int_0^L \phi_{\htilde^L}(z,\{ 0,y\})dz = \int_0^L \bigl( \htilde^L(0,z) + \htilde^L(y,z) - \htilde^L(0,z) \htilde^L(y,z) \bigr)dz. 
\end{equation}
Now, by the Schwarz inequality, 
\begin{equation} \label{eqn:cs_ineq}
\Bigl( \int_0^L \htilde^L(0,z) \htilde^L(y,z) dz \Bigr)^2 \leq 
\Bigl( \int_0^L (\htilde^L(0,z)^2) dz \Bigr)
\Bigl( \int_0^L (\htilde^L(y,z)^2) dz \Bigr).
\end{equation} 
As $\htilde^L(y,\cdot)$ is just a circular shift of $\htilde^L(0,\cdot)$, we also have 
\begin{equation} \label{eqn:circular}
\int_0^L \htilde^L(y,z)dz = \int_0^L \htilde^L(0,z)dz, \quad 
\int_0^L (\htilde^L(y,z)^2)dz = \int_0^L (\htilde^L(0,z)^2)dz.
\end{equation}
Substituting \eqref{eqn:cs_ineq} and \eqref{eqn:circular} into \eqref{eqn:phi_integral_0}, we get
\begin{equation} \label{eqn:phi_integral}
    \int_0^L \phi_{\htilde^L}(z,\{ 0,y\})dz  
    \geq 2\int_0^L \htilde^L(0,z)dz - \int_0^L (\htilde^L(0,z))^2 dz,
\end{equation}
Now, 
\begin{equation} \label{eqn:h1norm}
\begin{aligned}
\frac{1}{R_L} \int_0^L \htilde^L(0,z)dz 
&= \frac{2}{R_L} \int_0^{R_L^{1+1/\alpha}} H\Bigl( \frac{z}{R_L} \Bigr)dz \\
&= 2 \int_0^{R_L^{\alpha}} H(x)dx \xrightarrow{L\to \infty} 2\| H\|_1,
\end{aligned}
\end{equation}
since $R_L$ tends to infinity as $L$ does. A similar calculation yields 
\begin{equation} \label{eqn:h2norm}
\frac{1}{R_L} \int_0^L (\htilde^L(0,z))^2dz \xrightarrow{L\to \infty} 2R_L\| H\|_2^2.
\end{equation}
Substituting \eqref{eqn:h1norm} and \eqref{eqn:h2norm} in \eqref{eqn:phi_integral}, we get  
$$
\liminf_{L\to \infty} \frac{1}{R_L} \int_0^L \phi_{\htilde^L}(z,\{ 0,y\}dz \geq 4\| H\|_1 - 2\| H\|_2^2.
$$
Recall that, by assumption, there exists $\delta>0$ such that $\| H\|_2^2 \leq \| H\|_1-\delta$. Hence, it follows that 
$$
\liminf_{L\to \infty} \frac{1}{R_L} \int_0^L \phi_{\htilde^L}(z,\{ 0,y\} dz \geq 2\| H\|_1+2\delta.
$$
Since $R_L=\ln(\tau L)/2\| H\|_1$, we now obtain from \eqref{eqn:b2_bd} that 
\begin{align*}
\limsup_{L\to \infty} \limsup_{m\to \infty} b_2 &\leq \limsup_{L\to \infty} 6L R_L^{1+1/\alpha} e^{-\ln(\tau L) -2\delta R_L} \\ 
&= \limsup_{L\to \infty} \frac{6}{\tau} R_L^{1+1/\alpha} e^{-2\delta R_L} = 0,
\end{align*}
since $R_L$ tends to infinity with $L$.
\end{proof}

\begin{proof}[Proof of Theorem \ref{thm:IsoNodesPoisson_torus}]

Observe from Lemmas \ref{lem:Chen-Stein}, \ref{lem:b3_bound}, \ref{lem:b1_bound} and \ref{lem:b2_bound} that 
$$
\lim_{L\to \infty} \lim_{m\to \infty} d_{TV}(\iso, Po(\E[\iso]) = 0.
$$
It follows from Lemmas \ref{lem:niso_nisotilde_lim} and \ref{lem:discrete_approx}, and the triangle inequality, that 
$$
\lim_{L\to \infty} \lim_{m\to \infty} d_{TV}(\niso, Po(\E[\iso]) = 0.
$$
It is straightforward to check that $\E[\iso]=mLp_i$ tends to $1/\tau$. Finally, a straightforward calculation shows that, if a sequence $\lambda_n$ converges to $\lambda$, then $Po(\lambda_n)$ converges to $Po(\lambda)$ in total variation distance. Hence, invoking the triangle inequality once more, we obtain that 
$$
d_{TV}(\niso,Po(1/\tau)) \to 0 \mbox{ as } m,L\to \infty.
$$
This completes the proof of the theorem.
\end{proof}

{\bf Acknowledgements} MW was supported by the EPSRC Centre for Doctoral Training in Communications (EP/I028153/1 and EP/L016656/1).

\bibliographystyle{apalike}
\bibliography{bibliography.bib}

\end{document}